\documentclass[12pt, a4paper, reqno]{amsart} 

\usepackage{amssymb}
\usepackage{enumerate}

\newtheorem{theorem}{Theorem}[section]

\newtheorem{claim}[theorem]{Claim}
\newtheorem{problem}{Problem}
\newtheorem{alemma}[theorem]{Amalgamation Lemma}

\theoremstyle{definition}

\newtheorem{definition}[theorem]{Definition}

\newcommand{\mc}[1]{\mathcal{#1}}
\def\<{\left\langle}
\def\>{\right\rangle}
\def\br#1;#2;{\bigl[ {#1} \bigr]^ {#2} }

\newcommand{\force}{\Vdash}

\def\qedref#1{$\qed_{\ref{#1}}$}

\def\sbar{\underline{\sigma}}
\def\sstr{{\sigma}^*}

\author[S. Shelah]{Saharon Shelah}
\thanks{Publication 588} 
\address{Institute of Mathematics\\Hebrew University, Jerusalem }
 \email{shelah@math.huji.ac.il}
\title{Large weight does not yield an irreducible base}
\date{\today}
\begin{document}
\begin{abstract}
Answering a question of Juh\'asz, Soukup and Szentmiklóssy, we show that it is consistent that 
some  first countable space of uncountable weight 
does not contain an uncountable subspace which has an irreducible base.
\end{abstract}
\maketitle
\section{Introduction}

For a topological space $X, w(X)$ is the minimal cardinality of a base for $X$, $\chi (p, X) =\min\{ |u|: u$ is a neighbourhood base of $p  \}$, and $\chi (X) = \sup\{ \chi(p, X) : p \in X\}$.

In \cite{JSSz} the following problem was investigated: 
What makes a space have weight larger than its character?
The notion of {\em irreducible base} was introduced, and it was proved
\cite[Lemma 2.6]{JSSz} that {\em if a topological space $X$ has an irreducible base then
$w(X)=|X|\cdot \chi(X)$.}
 The  following question was formulated:
\begin{problem}
Does every  first countable space of uncountable weight 
contain an uncountable subspace which has  an irreducible base?
\end{problem}

We show that the answer is consistently NO. We thank Lajos Soukup for actually writing the paper.

\begin{definition}\label{df:irred}
Let $X$ be a topological space. 
A base ${\mc U}$ of $X$ is called {\em irreducible} 
if it has an {\em irreducible decomposition} ${\mc U}=\bigcup\{{\mc U}_x:x\in X\}$,
i.e,  $(i)$ and $(ii)$ below hold:
\begin{enumerate}[(i)]
\item ${\mc U}_x$ is a neighbourhood base of $x$ in $X$ for each $x\in X$.
\item for each $x\in X$   the family 
${\mc U}^-_x=\bigcup\limits_{y\ne x}{\mc U}_y$ is not a base of $X$.
\end{enumerate}
\end{definition}

\begin{theorem}
There is a c.c.c poset $P=\<P,\le\>$ of size ${\omega_1}$
such that in $V^P$ there is a   first countable space 
$X=\<{\omega_1},\tau\>$
of uncountable weight 
which does  not contain an uncountable subspace which has an irreducible base.
\end{theorem}

\noindent{\sl Proof.}
The elements of the poset $P$ will be finite 
``approximations'' of a base
$\{U({\alpha},n):{\alpha}<{\omega}_1,n<{\omega}\}$ of $X$.

We define the poset $P=\<P,\le\>$ as follows.
The underlying set of  $P$ consists of the 
 triples $\<A,n,U\>$ satisfying (P\ref{p:1})--(P\ref{p:3}) below:
\begin{enumerate}[(P1)]
\item \label{p:1} $A\in \br {\omega_1};<{\omega};$, $n\in {\omega}$ and
$U$ is a function, $U:A\times n\to \mc P(A)$,
\item \label{p:2}${\alpha}\in U({\alpha},i)\subset U({\alpha},i-1)$ for each 
${\alpha}\in A$ and $i<n$,
\item \label{p:3} 
If ${\beta}\in U({\alpha},i)\subset U({\beta},0)$  for some  $i<n$, then $\beta\le \alpha$.
\end{enumerate}
For $p\in P$  write $p=\<A_p,n_p,U_p\>$.
Let us remark that property (P3) will guarantee that $\mathrm{w}(X)=\omega_1$.

Define the order $\le$ on $P$ as follows.
For $p,q\in P $ we put $q\le p$ if
\begin{enumerate}[(a)]
\item $A_p\subset A_q$,
\item $n_p\le n_q$,
\item $U_p({\alpha},i)=U_q({\alpha},i)\cap A_p$ for each 
$\<{\alpha},i\>\in A_p\times n_p$,
\item\label{x}for each 
$\<{\alpha},i\>, \<{\beta},j\>\in A_p\times n_p$ ,
\begin{gather}
\tag{\ref{x}1 }  \text{ if $U_p({\alpha},i)\cap U_p({\beta},j)=\emptyset$ then 
$U_q({\alpha},i)\cap U_q({\beta},j)=\emptyset$,}  \\
\tag{\ref{x}2 } \text{if $U_p({\alpha},i)\subset U_p({\beta},j)$ then
$U_q({\alpha},i)\subset U_q({\beta},j)$.} 
\end{gather}
\end{enumerate}

 We say that the conditions $p_0=\<A_0,n_0,U_0\>$ and 
$p_1=\<A_1,n_1,U_1\>$
are {\em twins} iff  $n_0=n_1$, $|A_0|=|A_1|$ and 
denoting by ${\sigma}$ the unique $<_{\text{On}}$-preserving bijection between $A_0$ and
$A_1$ we have
\begin{enumerate}[({I}1)]
\item ${\sigma}\restriction {A_0\cap A_1}=\mathrm{id}_{A_0\cap A_1}$,
\item ${\sigma}$ is an isomorphism between $p_0$ and $p_1$, i.e.
for each ${\alpha}\in A_0$ and $ i<n_0$ we have
$U_1({\sigma}({\alpha}),i)={\sigma}''U_0({\alpha},i)$.
\end{enumerate}
We say that ${\sigma}$ is the {\em twin function} between $p_0$ and $p_1$.
Define the {\em smashing function} $\sbar$ of $p_0$ and $p_1$ as follows:
 $\sbar={\sigma}^{-1}\cup \mathrm{id}_{A_0}$.  
The function $\sstr$ defined by the formula 
$\sstr={\sigma}\cup {\sigma}^{-1}$ is called the
{\em exchange function} of $p_0$ and $p_1$.

The burden of the proof is to verify the  next lemma.

\begin{alemma}
\label{lm:twins}
Assume that  $p_0=\<A_0,n_0,U_0\>$ and   $p_1=\<A_1,n_1,U_1\>$ are twins,
$A_0\cap A_1< A_0\setminus A_1 < A_1\setminus A_0$, 
${\xi}_0\in A_{0}\setminus A_{1}$, ${\xi}_1={\sigma}({\xi}_0)$, where 
${\sigma}$ is the twin function between $p_0$ and $p_1$,   and let 
$k<m<n_0$. Then $p_0$
and $p_1$ have a common extension $p=\<A,n,U\>$ in $P$ such that 
\begin{equation}
\tag{$*$} {\xi}_0\in U({\xi}_1,m)\subset U({\xi}_1,k)\subset U({\xi}_0,k).
\end{equation}
\end{alemma}

\begin{proof}
Write $n=n_0=n_1$,
$D=A_0\cap A_1$ and ${{A^*}}=A_0\cup A_1$. 
Unfortunately we can not assume that $A={{A^*}}$ because in this case we
can not guarantee (P\ref{p:3}) for $p$. So we need to add further elements
to ${{A^*}}$ to get a large enough $A$ as follows. 
Choose a set $B\subset {\omega_1}\setminus {{A^*}}$ 
of cardinality $|{{A^*}}\times n|$ and fix a bijection 
$\rho$ between ${{A^*}}\times n$ and $B$. We will take $A={{A^*}}\cup B$. To simplify 
the notation we will write $\<{\alpha},i\>$ for $\rho({\alpha},i)$, 
for all ${\alpha}\in {{A^*}}$ and $i<n$, i.e. we identify  the elements of
$B$ and of ${{A^*}}\times n$.

The idea of the proof is the following: 
for each $\<{\alpha},i\>\in {{A^*}}\times n$ we put 
the element $\<{\alpha},i\>$ into $U({\alpha},i)$. On the other hand, 
we try to keep $U(\alpha, i)$ small, so
  we put  $\<\beta,j\>$ into $U({\alpha},i)$ if
and only if we can ``derive'' from the property (d2) that $U({\beta},j)\subset U(\alpha,i)$ should
hold in any condition $p=\<A,n,U\>$ which is a common extension of
$p_0$ and $p_1$ and which satisfies $(*)$.

The condition $p$ will be constructed in two steps.
First we construct a condition $p'=\<A,n,U'\>$ extending both $p_0$ and
$p_1$. This $p'$ can be considered as the minimal amalgamation of 
$p_0$ and $p_1$. Then, in the second step, we carry out  small
modifications on the function $U'$, namely we increase its value
on certain places to guarantee $(*)$.

Now we carry out our construction. 
For $\varepsilon<2$ and $\<{\beta},j\>\in A_\varepsilon\times n$ let
\begin{equation}\label{eq:v}
V_\varepsilon({\beta},j)=\left\{
\<{\alpha},i\>\in A_\varepsilon\times n: U_\varepsilon({\alpha},i)\subset U_\varepsilon({\beta},j)
\right\}
\end{equation}
and
\begin{equation}\label{eq:w}
\begin{split}
W_\varepsilon({\beta},j)=\{
\<{\alpha},i\>\in &A_{1-\varepsilon}\times n :\ \exists
\<{\gamma},l\>\in D\times n\\& 
U_{1-\varepsilon}({\alpha},i)\subset U_{1-\varepsilon}({\gamma},l)\land
U_\varepsilon({\gamma},l)\subset U_\varepsilon({\beta},j)
\}
\end{split}
\end{equation}
If we want to define $p'$ in such  a way that $p'\le p_0, p_1$, then (d2) implies that 
$U'(\alpha,i)\subset U'(\beta,j)$ should hold whenever $\<\alpha,i\>\in V(\beta,j)\cup W(\beta,j)$.

Now we are ready to define the function $U'$.
For $\varepsilon<2$, ${\beta}\in A_\varepsilon$ and $j<n$ let
\begin{equation}\label{eq:u'_old}
U'({\beta},j)=U_\varepsilon({\beta},j)\cup U_{1-\varepsilon}(\sstr({\beta}),j)\cup V_\varepsilon({\beta},j)
\cup W_\varepsilon({\beta},j).
\end{equation}
For $\<{\alpha},i\>\in {{A^*}}\times n$ and $j<n$ let
\begin{equation}\label{eq:u'_new}
U'(\<{\alpha},i\>,j)=\{\<{\alpha},i\>\}.
\end{equation}
Let us remark that $U'(\delta,j)$ is well-defined even for $\delta\in A_0\cap A_1$. Indeed, 
in this case  $\sstr(\delta)=\delta$ and $V_\varepsilon(\delta,j)=W_{1-\varepsilon}(\delta,j)$,
and so 
\begin{equation}\notag
U'({\delta},j)=U_0({\delta},j)\cup U_{1}({\delta},j)
\cup V_0({\delta},j)\cup V_1({\delta},j). 
\end{equation} 
Now put
\begin{equation}\notag
p'=\<A,n,U'\>.
\end{equation}

\begin{claim}\label{cl:push}
If ${\alpha}\in U'({\beta},j)$ then $\sbar(\alpha)\in U_0(\sbar(\beta),j)$.
\end{claim}
Indeed, if ${\beta}\in A_\varepsilon$ then 
$U'({\beta},j)\cap {A^*} =U_\varepsilon({\beta},j)\cup U_{1-\varepsilon}(\sstr({\beta}),j)$.

\begin{claim}\label{cl:push2}
If $\<{\alpha},i\>\in U'({\beta},j)$ then 
$\sbar(\alpha)\in U_0(\sbar(\beta),j)$.
\end{claim}
\noindent{\sl Proof of the  Claim}.
Assume that ${\beta}\in A_\varepsilon$.
If $\<\alpha,i\>\in V_\varepsilon(\beta,j)$ then 
 $\alpha\in U_\varepsilon({\alpha},i)\subset U_\varepsilon({\beta},j)$
and  $ U_\varepsilon({\beta},j)\subset U'(\beta,j)$.
So we have ${\alpha}\in U'({\beta},j)$ which implies 
$\sbar(\alpha)\in U_0(\sbar(\beta),j)$ by Claim \ref{cl:push}.

If  $\<\alpha,i\>\in W_{1-\varepsilon}(\beta,j)$ then  
for some
$\<{\gamma},l\>\in D\times n$ we have  
$U_{1-\varepsilon}({\alpha},i)\subset U_{1-\varepsilon}({\gamma},l)\land
U_\varepsilon({\gamma},l)\subset U_\varepsilon({\beta},j)$.
Thus  $\alpha\in  U_\varepsilon({\beta},j)\subset U'(\beta,j)$,
which implies 
$\sbar(\alpha)\in U_0(\sbar(\beta),j)$ by Claim \ref{cl:push}.
\hfill \qedref{cl:push2}

\begin{claim}\label{cl:p'1}
$p'\in P$.
\end{claim}

\noindent{\sl Proof of the claim \ref{cl:p'1}. } 
(P\ref{p:1}) and (P\ref{p:2}) clearly hold, so we need to check only (P\ref{p:3}). 

Assume on the contrary that (P\ref{p:3}) fails for $p'$.
Since $U'(\<\nu,s\>,j)=\{\<{\nu},s\>\}$ 
by (\ref{eq:u'_new})
for each $\<{\nu},s\>\in B$  and $j<n$, we can assume that some 
${\alpha}<{\beta}\in {{A^*}}$ and 
$i<n$ witness that (P\ref{p:3}) fails, i.e. ${\beta}\in U'({\alpha},i)\subset U'({\beta},0)$.
Then $\sbar(\beta)\in U_0(\sbar(\alpha),i)\subset U'(\sbar(\beta),0)$ by Claim \ref{cl:push}.
Since $p_0$ satisfies (P\ref{p:3}) it follows that  $\sbar(\beta)\le \sbar(\alpha)$, and so
${\alpha}\in A_0\setminus A_1$ and ${\beta}\in A_1\setminus A_0$.
Consider the element $u=\<{\alpha},i\>\in A\setminus {{A^*}}$. 
Then $u\in U'({\alpha},i)$ and
so $u\in U'({\beta},0)$ as well. By the definition of $U'({\beta},0)$ 
this means  that $\<\alpha,i\>\in W_1(\beta,0)$, that is,   
there is $\<{\gamma},l\>\in D\times n$
such that $U_0({\alpha},i)\subset U_0({\gamma},l)$ and 
$U_1({\gamma},l)\subset U_1({\beta},j)$. Thus 
\begin{equation}
\sbar(\beta)\in U_0({\alpha},i)\subset U_0({\gamma},l)\subset U_0(\sbar(\beta),0)
\end{equation} 
by Claim \ref{cl:push}.
Thus  $\sbar(\beta)\in U_0({\gamma},l)\subset U_0(\sbar(\beta),0)$
 and so $\sbar(\beta)\le \gamma $ because $p_0$ satisfies (P\ref{p:3}).
But this is a contradiction because  $\gamma\in D=A_0\cap A_1$, 
$\sbar(\beta)\in A_0\setminus A_1$ and we assumed that $(A_0\cap A_1)<(A_0\setminus A_1)$.

\hfill\qedref{cl:p'1}

%
%
%
%

\begin{claim}\label{cl:p'}
$p'\le p_0,p_1$.
\end{claim}

\noindent{\sl Proof of  claim \ref{cl:p'}. }
Conditions (a) and (b) are clear.

To check (c) assume that $\alpha\in A_\varepsilon$ and $i\in  n$.
By (\ref{eq:u'_old}),
\begin{multline}\notag
 U'(\alpha,i)\cap A_\varepsilon=  
(U_\varepsilon (\alpha,i)\cup U_{1-\varepsilon }(\alpha, i))\cap A_\varepsilon  =\\
U_\varepsilon (\alpha,i)\cup ( U_{1-\varepsilon }(\alpha, i) \cap A_\varepsilon)=
U_\varepsilon (\alpha,i)
\end{multline}
because $U_{1-\varepsilon }(\alpha, i)=\sstr[U_{\varepsilon }(\alpha, i)]$.

To check (d1) assume that $\beta,\gamma\in A_\varepsilon$ and  $j,k<n$ 
such that 
$U'(\beta,j)\cap U'(\gamma,k)\ne \emptyset$. Fix $x\in U'(\beta,j)\cap U'(\gamma,k)$. 
Then 
\begin{equation}\notag
 \sbar(\alpha)\in U_0(\sbar(\beta),j)\cap U_0(\sbar(\gamma),k)
\end{equation} 
by  Claim \ref{cl:push} if $x=\alpha\in {A^*}$, and 
by Claim \ref{cl:push2}
if $x=\<\alpha,i\>\in A\setminus {A^*}$.

If $\varepsilon=0$ then $\sbar(\beta)=\beta$ and $\sbar(\gamma)=\gamma$, so 
$\sbar(\alpha)\in U_\varepsilon(\beta,j)\cap U_\varepsilon(\gamma,k)$.

If $\varepsilon=1$ then $\sbar(\beta)=\sstr(\beta)$ and $\sbar(\gamma)=\sstr(\gamma)$, 
and so $\sstr(\sbar(\alpha))\in U_\varepsilon(\beta,j)\cap U_\varepsilon(\gamma,k)$.

Finally to check (d2) assume that $\beta,\gamma\in A_\varepsilon$ and  $j,k<n$ 
such that $U_\varepsilon(\beta,j)\subset  U_\varepsilon(\gamma,k)$. 
Then clearly 
\begin{equation}\notag
 U_{1-\varepsilon}(\beta,j)=\sstr[U_\varepsilon(\beta,j)]\subset
\sstr[U_\varepsilon(\gamma,k)]=
  U_{1-\varepsilon}(\gamma,k),
\end{equation} 
moreover 
$V_\varepsilon(\beta,j)\subset  V_\varepsilon(\gamma,k)$ by (\ref{eq:v}), and
$W_\varepsilon(\beta,j)\subset  W_\varepsilon(\gamma,k)$ by (\ref{eq:w}), and so
$U'(\beta,j)\subset  U'\varepsilon(\gamma,k)$ by  (\ref{eq:u'_new}). 
\hfill\qedref{cl:p'}


Now carry out the promised modification of $U'$ to obtain  $U$ as follows.
If $z\in A$ and $j<n$ let
\begin{equation}\notag
U(z,j)=\left\{
\begin{array}{ll}
U'(z,j)\cup U'({\xi}_1,k)&\text{if 
$U_0({\xi}_0,k)\subset U_0(z,j)$,}\\
U'(z,j)&\text{otherwise.}
\end{array}
\right.
\end{equation}

Put
\begin{equation}\notag
p =\<A,n,U \>.
\end{equation}

If $U_0({\xi}_0,k)\subset U_0(z,j)$
then $U_1({\xi}_1,k)\subset U_1(\sstr (z),j)\subset U'(z,j)$
and $W_1(\xi_1, k)\subset V_0(\xi_0,k)\subset U'(z,j)$.
So 
\begin{equation}\label{eq:u-u'}
U(z,j)\setminus U'(z,j)\subset V_1(\xi,k). 
\end{equation} 
Moreover  
\begin{equation}\label{eq:u2}
U(z,j)=\left\{
\begin{array}{ll}
U'(z,j)\cup V_1({\xi}_1,k)&\text{if 
$U_0({\xi}_0,k)\subset U_0(z,j)$,}\\
U'(z,j)&\text{otherwise.}
\end{array}
\right.
\end{equation}

\begin{claim}\label{cl:push3}
If $\<{\alpha},i\>\in U({\beta},j)$ then 
$\sbar(\alpha)\in U_0(\sbar(\beta),j)$.
\end{claim}

Indeed, if $\<{\alpha},i\>\in U({\beta},j)$ then 
$\<{\alpha},i\>\in U'({\beta},j)$ or $\<{\alpha},i\>\in U'(\sstr({\beta}),j)$, 
and now apply Claim \ref{cl:push2}.

\begin{claim}\label{cl:p}
$p\in P$.
\end{claim}

\noindent{\sl Proof of claim \ref{cl:p}. }
(P\ref{p:1}) and (P\ref{p:2}) clearly hold, so we need  to check (P\ref{p:3}) only.
 
Assume on the contrary that (P\ref{p:3}) fails for $p$.
Since $U(\<{\nu},s\>,j)=\{\<{\nu},s\>\}$ for each $\<{\nu},s\>\in A\setminus A^* $ and
$j<n$ we can assume that there are ${\alpha}<{\beta}\in {{A^*}}$ and 
$i<n$ witness that (P\ref{p:3}) fails, i.e.  
\begin{equation}\label{eq:ind}
{\beta}\in U({\alpha},i)\subset U({\beta},0).
\end{equation} 
Then $\sbar(\beta)\in U_0(\sbar(\alpha),i)\subset U(\sbar(\beta),0)$.
But $p_0$ satisfies (P\ref{p:3}) so $\sbar(\beta)\le \sbar(\alpha)$, and so
${\alpha}\in A_0\setminus A_1$ and ${\beta}\in A_1\setminus A_0$.
Thus $U_0(\beta,j)$ is undefined, and so 
\begin{equation}
U'({\beta},0)= U({\beta},0)  
\text{ and }  U({\alpha},i)\setminus U'({\alpha},i)\subset A\setminus A^*.
\end{equation} 
by (\ref{eq:u2}). So (\ref{eq:ind}) yields 
\begin{equation}\notag
{\beta}\in U'({\alpha},i)\subset U'({\beta},0),
\end{equation}
However this is a contradiction 
because $p'$ satisfies (P\ref{p:3}).
\hfill\qedref{cl:p}


\begin{claim}\label{cl:p2}
$p\le p_0, p_1$.
\end{claim}

\noindent{\sl Proof. } 
(a) and (b) are trivial.
(c) also holds because $p'\le p_\varepsilon$ and
$(U(\alpha,i)\setminus U'(\alpha,i))\cap A_\varepsilon=\emptyset$ by (\ref{eq:u-u'}.)

To check (d1) assume that $\beta,\gamma\in A_\varepsilon$ and  $j,k<n$ 
such that $U(\beta,j)\cap U(\gamma,k)\ne \emptyset$. Pick
$x\in U(\beta,j)\cap U(\gamma,k)$. 
Then 
\begin{equation}\notag
 \sbar(\alpha)\in U_0(\sbar(\beta),j)\cap U_0(\sbar(\gamma),k)
\end{equation} 
by  Claim \ref{cl:push} if $x=\alpha\in {A^*}$, and 
by Claim \ref{cl:push3}
if $x=\<\alpha,i\>\in A\setminus {A^*}$.

If $\varepsilon=0$ then $\sbar(\beta)=\beta$ and $\sbar(\gamma)=\gamma$, so 
$\sbar(\alpha)\in U_\varepsilon(\beta,j)\cap U_\varepsilon(\gamma,k)$.

If $\varepsilon=1$ then $\sbar(\beta)=\sstr(\beta)$ and $\sbar(\gamma)=\sstr(\gamma)$, 
and so $\sstr(\sbar(\alpha))\in U_\varepsilon(\beta,j)\cap U_\varepsilon(\gamma,k)$.

Finally to check (d2) assume that $\beta,\gamma\in A_\varepsilon$ and  $i,j<n$ 
such that $U_\varepsilon(\beta,i)\subset  U_\varepsilon(\gamma,j)$. 
Since $p'\le p_\varepsilon$ we have  $U'(\beta,i)\subset  U'(\gamma,j)$.
If $U(\beta,i)=U'(\beta,i)$, we are done.
So we can assume that $U(\beta,i)=U'(\beta,i)\cup V(\xi_1,k)$.
Then $\varepsilon=0$ and $U_0(\xi_0,k)\subset U_0(\beta,i)$.
But then $U_0(\xi_0.k)\subset U_0(\gamma,j)$ and so 
$U(\gamma,j)=U'(\gamma,j)\cup V(\xi_1,k)$, and so 
$U(\beta,i)\subset  U(\gamma,j)$.
\hfill\qedref{cl:p2}

Since $p$ satisfies $(*)$, the amalgamation lemma is proved.
\hfill\qedref{lm:twins}

Using the amalgamation lemma it is easy to complete the proof of the theorem.

By standard $\Delta$-system argument, any uncountable set of conditions contains
two elements, $p_0$ and $p_1$, which are twins. So, by Lemma \ref{lm:twins},
they have a common extension $p$. So $P$ satisfies c.c.c.

If $\mc G$ is a generic filter,
for $\alpha<{\omega_1}$ and $i<\omega $ put
\begin{equation}
U(\alpha,i)=\cup\{U_p(\alpha,i):p\in \mc G, \alpha\in A_p, i<n_p\},
\end{equation} 
and let $\mc U_\alpha=\{U(\alpha,i):i<\omega\}$ be the base of the point $\alpha$
in $X=\<{\omega_1},\tau\>$.

By (P\ref{p:3}), a countable 
subfamily of $\{U(\alpha,i):\alpha<{\omega_1}, i<\omega\}$
is not a base of $X$. So $\mathrm{w}(X)={\omega}_1$.

Finally we show that $X$
does  not contain an uncountable subspace which has  an irreducible base.

Assume on the contrary that
\begin{multline}
\notag
r\force\text{the subspace $\dot Y=\{\dot y_\xi:\xi<{\omega_1}\}$  has an irreducible base $\mc B$,}
\\
\text{and $\{\dot {\mc B}_{y_{\xi}}:{\xi}<{\omega}_1\}$ is an 
irreducible decomposition of $\dot{\mc B}$}.
\end{multline}
We can assume that $r\force\dot y_\xi\ge \check\xi$.

For each ${\xi}<{\omega}_1$ pick a condition $r_{\xi}$ and 
$k_{\xi}\in {\omega}$ such that 
\begin{equation}\label{eq:irr}
\text{
$r_{\xi}\force$ ``if 
$V\in \mc B$ with  $\dot y_{\xi}\in V\subset U(\dot y_{\xi},\check k_{\xi})$ then 
$V\in \mc B_{y_\xi}$''.
} 
\end{equation}

For each ${\xi}<{\omega}_1$ pick a condition $p_{\xi}\le r_{\xi}$, an ordinal
$\alpha_\xi\ge \xi $, 
a name $\dot V_{\xi}$ and a natural number $m_{\xi}<{\omega}$ such that
$\alpha_\xi\in A_{p_\xi}$ and   
\begin{equation}\label{eq:between}
p_\xi\force
\dot y_\xi=\check \alpha_\xi, 
\dot V_{\xi}\in \dot{\mc B}_{\alpha_{\xi}} \text{ and } 
U(\check \alpha_{\xi},\check m_{\xi})\subset  
\dot V_{\xi}\subset U(\check \alpha_{\xi},\check k_\xi). 
\end{equation}

By standard argument find $I\in \br {\omega_1};{\omega_1};$ such that 
\begin{enumerate}[(i)]
 \item $m_\xi=m$ and $k_\xi=k$ for each $\xi\in I$,
\item the sequence $\{\alpha_\xi:\xi\in I\}$ is strictly increasing,
\item the conditions $\{p_\xi:\xi\in I\}$ are pairwise twins,
\item  $\sigma_{\xi,\eta}(\alpha_\xi)=\alpha_\eta$ for
$\{\xi,\eta\} \in \br I;2;$, where $\sigma_{\xi,\eta}$ is the twin function.
\end{enumerate}

Pick $\xi<\eta$ from $I$.
By the Amalgamation Lemma there is a common extension $p$ of $p_\xi$ and $p_\eta$
such that 
\begin{equation}
 \alpha_\xi\in  U_p(\alpha_\eta,m) \land U_p(\alpha_\eta,k)\subset U_p(\alpha_\xi,k). 
\end{equation} 
Then, by (d2),
\begin{equation}
p\force \check \alpha_\xi\in  U(\check \alpha_\eta,\check m) \land 
U(\check \alpha_\eta,\check k)\subset U(\check \alpha_\xi,\check k).
\end{equation} 
Then, by (\ref{eq:between}),
\begin{equation}
p\force \dot V_\eta\in \mc B_{\alpha_\eta}\text{ and }
\check \alpha_\xi\in  U(\check \alpha_\eta,\check m)\subset
\dot V_\eta\subset    
U(\check \alpha_\eta,\check k)\subset U(\check \alpha_\xi,\check k),
\end{equation} 
which contradicts (\ref{eq:irr}).

This completes the proof of the Theorem.
\end{proof}

\end{document}